\documentclass[12pt]{article}

\usepackage{graphics,graphicx,amssymb,amsthm,amsmath,latexsym}

\begin{document}

\title{Numerical and analytical investigation of the free boundary
confluence for the phase field system.}

\author{V.G.~Danilov\thanks{Moscow Technical University of
Communication and Informatics,\hfill\break e-mail: danilov@miem.edu.ru },\ \ \
V.Yu.~Rudnev\thanks{Moscow Technical University of Communication
and Informatics,\hfill\break e-mail: vrudnev78@mail.ru }}
\date{}
\maketitle

\begin{abstract}
In this paper we numerically research the solutions of the phase field system for the
spherically symmetric Stefan-Gibbs-Thomson problem in the case of interaction of the free
boundaries. We analyze the effect of the soliton type disturbance of the temperature
in the point of the contact of the free boundaries.
\end{abstract}

\section{INTRODUCTION. STATEMENT OF THE PROBLEM}

The main goal of this paper is the numerical research of the effect of the soliton type disturbance
of the temperature in the point of the contact of the free boundaries. This effect we consider in
the case of the phase field model for the Stefan-Gibbs-Thomson problem. The difference between this
problem and the classical Stefan problem is that the surface tension is taken into account in the
Stefan-Gibbs-Thomson problem. At the beginning of the paper we briefly give the main analytic
results about the confluence of the free boundaries and then we illustrate these results by computer
simulation.

\begin{figure}[ht!]
\begin{center}
\includegraphics{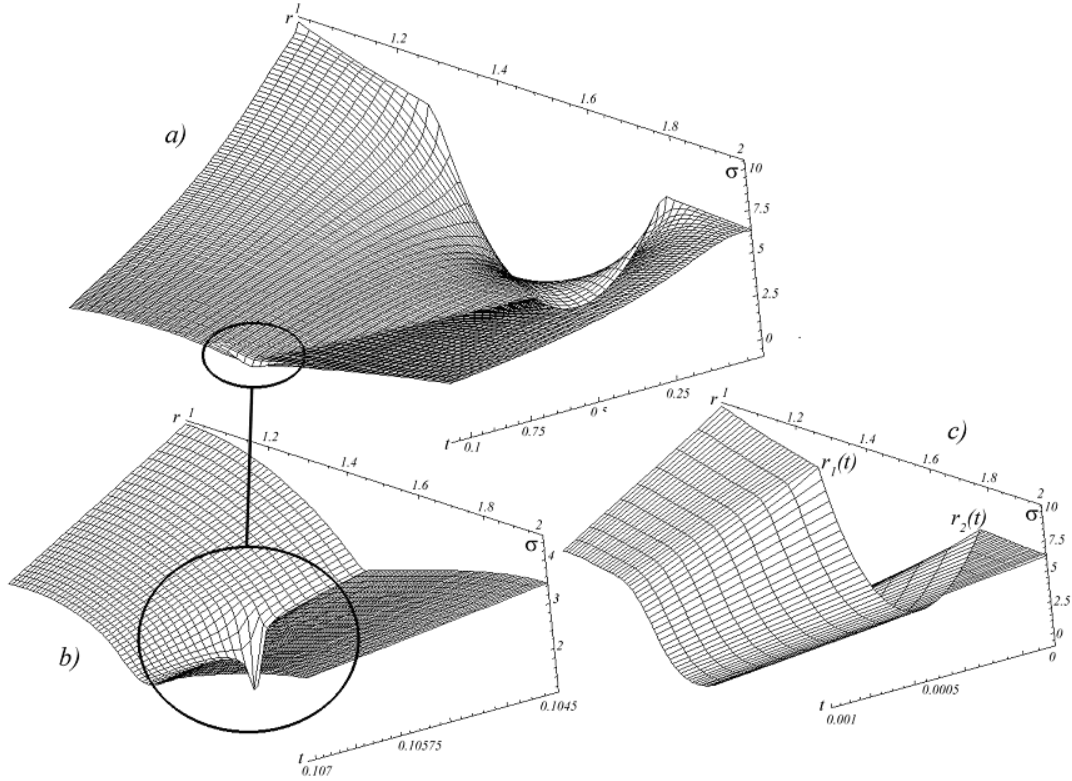}
\end{center}
\caption{Effect of the soliton type disturbance of the temperature
($\varepsilon=0.003$). a) Coordinate and time dependence of the
temperature. b) Soliton type disturbance of the temperature in the
neighborhood of the point of the contact of the free boundaries.
c) Initial data of the temperature and the evolution of the free
boundaries $r_{1}(t)$, $r_{2}(t)$ in prime.}\label{F0}
\end{figure}

The effect of the soliton type (negative) disturbance of the
temperature $\sigma$ of the free boundaries in the point of the
contact is shown in Fig. \ref{F0}.

\begin{figure}
\begin{center}
\includegraphics[width=7cm]{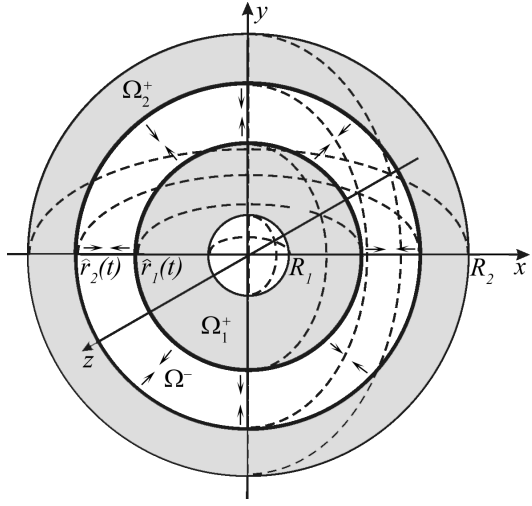}
\end{center}
\caption{The free boundaries in the spherically symmetric case.}\label{F1}
\end{figure}

We present the results of the computer simulation for the phase field system in the spherically
symmetric case. Namely, we consider the domain $Q=\Omega\times[0,t_{1}]$, where
$\Omega=[R_{1},R_{2}]$ is the spherical layer in the spherical coordinates. We assume that the
domain $\Omega$ is divided into the three layers $\Omega^{+}_{1,2}(t)$ and $\Omega^{-}(t)$ as
follows:
\begin{equation*}
\Omega^{+}_{1}(t)=[R_{1},r_{1}(t)], \quad
\Omega^{-}(t)=[r_{1}(t),r_{2}(t)], \quad \Omega^{+}_{2}(t)=[r_{2}(t), R_{2}],
\end{equation*}
where $r_{i}(t)=\Gamma_{i}(t)$, $i=1,2$ is the free boundaries between phases "$+$"\ and "$-$". We
assume that the phase~"$+$"\ occupies the layer $\Omega^{+}_{1,2}(t)$, and the phase~"$-$"\ occupies
the layer $\Omega^{-}(t)$ (see Fig.~\ref{F1}).

In this case the phase field system  have the form \cite{Cag1}
\begin{equation}\label{PF:1}
L\theta= -\frac{\partial u}{\partial t},\qquad
\varepsilon Lu-\frac{u-u^{3}}{\varepsilon}-\varkappa\theta=0,
\end{equation}
where $\varkappa=\sqrt{2}/3$,
\begin{equation*}\label{L:1}
L=\frac{\partial}{\partial t}-\frac{1}{r^{2}}\frac{\partial}{\partial r}
\left(r^{2}\frac{\partial}{\partial r}\right),
\qquad r\in[R_{1},R_{2}], \quad t\in[0,t_{1}].
\end{equation*}

The function $\theta=\theta(r,t,\varepsilon)$ has the meaning of the temperature, and the function
$u=u(r,t,\varepsilon)$ (so-called the order function) determines the phase state of the medium.
Namely, the value $u\simeq-1$ corresponds to the phase~"$-$"\ in the layer $\Omega^{-}(t)$, and the
value $u\simeq 1$ corresponds to the the phase~"$+$"\ in the layers $\Omega^{-}_{1,2}(t)$.

We denote
\begin{equation}\label{sigma}
\theta(r,t)=\frac{\sigma(r,t)}{r}.
\end{equation}
The phase field system (\ref{L:1}) we can write in the form
\begin{align}
&\frac{\partial\sigma}{\partial
t}-\frac{\partial^{2}\sigma}{\partial r^{2}}=
-r\frac{\partial u}{\partial t},\label{PhHeat}\\
&\varepsilon
Lu=\frac{u-u^{3}}{\varepsilon}+\varkappa\frac{\sigma}{r}.\label{PhAC}
\end{align}

Passing to the limit as $\varepsilon\to 0$ in (\ref{PhHeat}), (\ref{PhAC}) we obtain the
Stefan-Gibbs-Thomson problem for the each free boundary (see \cite{Cag1, DOR, PS})
\begin{align}
&\frac{\partial\overline{\sigma}}{\partial t}
=\frac{\partial^{2}\overline{\sigma}}{\partial r^{2}},\qquad
r\in[R_{1},R_{2}],\quad r\neq \hat{r}_{i}(t),\quad
i=1,2,\label{SH:1}\\
&\left.\frac{\overline{\sigma}}{r}\right|_{r=\hat{r}_{i}(t)}=(-1)^{i+1}\left(
\hat{r}_{i}'(t)+\frac{2}{\hat{r}_{i}(t)}\right),\label{SGT:1}\\
&\left.\left[\frac{\partial\overline{\sigma}}{\partial r}\right]
\right|_{r=\hat{r}_{i}(t)}=(-1)^{i+1}2\hat{r}_{i}(t)\hat{r}'_{i}(t).\label{SSt:1}
\end{align}
This passage to the limit is possible, for example, in the case where the corresponding limit
problems have classical solutions. In this case, the weak limits as $\varepsilon\to 0$ of solutions
(\ref{PhHeat}), (\ref{PhAC}) give these solutions~\cite{Cag1, DOR, PS}. The existence of the
classical solution $\overline{\sigma}$ of problem (\ref{SH:1})--(\ref{SSt:1}) is discussed in
\cite{M}.

By $t^{*}\in(0,t_{1})$ we denote the instant of time of the confluence of the free boundaries, and
$r^{*}=\hat{r}_{1}(t^{*})=\hat{r}_{2}(t^{*})$ is the sphere of the contact of the free boundaries.

The smooth approximations (i.e. the approximate solution of the system (\ref{PhHeat}), (\ref{PhAC})
with misalignment that is small in the weak sense as $\varepsilon\to 0$, see \cite{DOR}) of the
Stefan-Gibbs-Thomson problem (of the phase field system) is constructed in \cite{D} in the one
dimensional case and this smooth approximations exists in the time included the instant of the
contact $t^{*}$. Moreover, this smooth approximations is constructed on the assumption of the
existence of the classical solution of the limit problem as $t\leqslant t^{*}-\delta$, where $\delta
> 0$ is any number. This approximations admits the passage to the limit as $\varepsilon\to 0$.

It is possible to show \cite{DOR} that in the common sense the asymptotic solution of the system
(\ref{PhHeat}),~(\ref{PhAC}) has the form
\begin{align}\label{AS:1}
\sigma^{as}_{\varepsilon}&=\bar{\sigma}^{-}(r,t)
+ \left(\bar{\sigma}^{+}(r,t)-\bar{\sigma}^{-}(r,t)\right)
\omega_{1}\left(\frac{\hat{r}_{2}(t)-r}{\varepsilon}\right)
\omega_{1}\left(\frac{r-\hat{r}_{1}(t)}{\varepsilon}\right),
\\
\label{AS:2}
u^{as}_{\varepsilon}&= 1+\omega_{0}\left(\frac{\hat{r}_{1}(t)-r}{\varepsilon}\right)
\\
&\qquad +\omega_{0}\left(\frac{r-\hat{r}_{2}(t)}{\varepsilon}\right)
+\varepsilon\left[\frac{\sigma^{as}_{\varepsilon}}{2}
+\omega\left(t,\frac{r-\hat{r}_{1}(t)}{\varepsilon},
\frac{r-\hat{r}_{2}(t)}{\varepsilon}\right)\right].\notag
\end{align}
as $t\leqslant t^{*}-\delta$, $\delta> 0$. Here $\omega_{1}(z)\to 0,1$ as $z\to\mp\infty$,
$\omega^{(k)}_{1}(z)\in\mathbb{S}(\mathbb{R}^{1}_{z})$ as $k>0$, $\hat{r}_{i}(t)$, $i=1,2$ is the
smooth functions, $r_{1}\leqslant r_{2}$, $\omega_{0}(z)=\tanh(z)$, and $\omega(t,z_{1},z_{2})\in
C^{\infty}([0,t^{*}];\,\mathbb{S}(\mathbb{R}^{2}_{z}))$. By $\mathbb{S}(\mathbb{R}^n)$ we denote the
Schwartz space of smooth rapidly decreases functions. If the initial data for (\ref{PhHeat}),
(\ref{PhAC}) has the form (\ref{AS:1}), (\ref{AS:2}) at $t=0$, then, for $t\leqslant t^{*}-\delta$,
we have the estimate
\begin{equation*}
\|u-u^{as}_{\varepsilon};C(0,T;L^{2}(\mathbb{R}^{1})\|
+\|\sigma- \sigma^{as}_{\varepsilon};\mathcal{L}^{2}(Q)\|\leqslant
c\varepsilon^{\mu},\qquad \mu\geqslant  3/2,
\end{equation*}
where $(\sigma,u)$ is a solution of system (\ref{PhHeat}),(\ref{PhAC}) (see~\cite{DOR} and
references in). Here $Q=\Omega\times[0,t^{*}-\delta]$, and the constant $c$ is independent
of~$\varepsilon$.

The main obstacle to the construction of approximations of solution in the case of confluence of
free boundaries is the fact that, instead of an ordinary differential equation whose solution is the
function $\omega_{0}(z)$~\cite{Cag1, DOR}, in the case of confluence of free boundaries, we must
deal with a partial differential equation for which the explicit form of the exact solution is
unknown.

In papers \cite{D, DR} using the the weak asymptotic method the solution of the phase field system
was constructed that describes the confluence of the free boundaries. Namely, the ansatz of the
order function has the form.
\begin{align}\label{AU:1}
 \check{u}&=\frac{1}{2}\left[1+
 \omega_{0}\left(\beta\frac{r_{1}^2-r^2}{\varepsilon}\right)+
 \omega_{0}\left(\beta\frac{r^2-r_{2}^2}{\varepsilon}\right)\right.\\
          &   -\left.
  \omega_{0}\left(\beta\frac{r_{1}^2-r^2}{\varepsilon}\right)
  \omega_{0}\left(\beta\frac{r^2-r_{2}^2}{\varepsilon}\right)\right]\notag.
\end{align}

The temperature is sought in the form
\begin{equation}\label{new23}
\check{\sigma}=e(r)\check{T}+q,
\end{equation}
$e(r)\in C_{0}^{\infty}([R_{1},R_{2}])$, $e\equiv 1$ for $r\in[\hat{r}_{1}(0),\hat{r}_{2}(0)]$. Here
$\check{T}$ is the model of the temperature, i.e. it is the function of the simplest structure which
describes the behaviors of the temperature qualitatively correct, and $q$ is an unknown smooth
function. Namely,
\begin{align}\label{ModT:1}
\check{T}&=\gamma_{1}^{+}(t)(r_{1}-r)H(r_{1}-r)
+\gamma_{2}^{+}(t)(r-r_{2})H(r-r_{2})\notag\\
&+\gamma^{-}(r,t)\frac{(r_{1}-r)(r-r_{2})}{r_{2}-r_{1}}
H(r-r_{1})H(r_{2}-r)\\
&+\hat{\gamma}(r,t)\frac{(r_{1}-r)(r-r_{2})}{r_{2}-r_{1}}
H(r_{1}-r)H(r-r_{2})+I(r,t).\notag
\end{align}
Here
\begin{align}
&\gamma^{-}=\frac{\gamma_{1}^{+}-\gamma_{2}^{-}}{2}
-(r^2-r^{*2})\frac{\gamma_{1}^{+}+\gamma_{2}^{-}}{\psi}\notag
\\
&\hat{\gamma}=\frac{\hat{\gamma}_{1}+\hat{\gamma}_{2}}{2}
-(r^2-r^{*2})\frac{\hat{\gamma}_{1}-\hat{\gamma}_{2}}{\psi}, \qquad
r^{*}=\frac{\sqrt{r_{1}^2+r_{2}^2}}{2}.\notag
\\
&I=\frac{k_{1}-k_{2}}{2} -(r^2-r^{*2})\frac{k_{1}+k_{2}}{\psi},
\quad k_{1}=r_{1}r_{1t}+2,\quad k_{2}=r_{2}r_{2t}+2,\label{I}
\end{align}
\begin{equation}\label{R1PH1} \psi(t,\varepsilon)=r_{2}^2(t,\varepsilon)-r_{1}^2(t,\varepsilon),
\end{equation}
So, from the given above formulas we see that the model $\check{T}$ is linear in $r$ in the layers
$\Omega^{+}_{i}(t)$, $i=1,2$ and parabolic in $r$ in the layer $\Omega^{-}(t)$.

In paper \cite{D} the formulas are obtained those determine the
functions contained in ansatzes (\ref{AU:1}), (\ref{new23}) in the
one dimensional case.

The analysis of these formulas shows that the contained in the phase field system temperature
$\theta$ has the soliton type disturbance in the neighborhood of the point of the contact. This
disturbance is localized in the space coordinates and in time. The "width" of this localization is
proportional to $\varepsilon$. In \cite{D} the limit as $\varepsilon\to 0$ of the amplitude of the
disturbance is derived. In the one dimensional case the value of this amplitude is given by formula
\begin{equation*}
\left.\left[\bar{\theta}\right]\right|_{x=x^*,\ t=t^*}
=-\lim_{t\to t^{*}-0}\lim_{\varepsilon\to 0}\frac{V_{1}+V_{2}}{2},
\end{equation*}
where $V_{1}$ and $V_{2}$ are the velocities of the merging
(one-dimensional) free boundaries.

The analytic treatment based on the weak asymptotics method shows
that such effect is general. For example, the exactly same formula
is correct in the case of merging globe layers \cite{DR} those is
shown in Fig.\ref{F1}. Namely,
\begin{equation*}
\left.\left[\bar{\theta}\right]\right|_{r=r^*,\ t=t^*}
=-\lim_{t\to t^{*}-0}\lim_{\varepsilon\to
0}\frac{V_{1n}+V_{2n}}{2},
\end{equation*}
$V_{1n}$ and $V_{2n}$ are the normal velocities of the merging three-dimensional free boundaries
(correspondingly, $\Gamma_{1}(t)$ and $\Gamma_{2}(t)$). If the free boundaries are asymmetrical,
then their principle curvatures have the like signs at the point of the contact. In this case the
amplitude of the jump of the temperature is determined by formula
\begin{equation*}
\left.\left[\bar{\theta}\right]\right|_{(x,y,z)=(x^*,y^*,z^*),\ t=t^*} =-\lim_{t\to
t^{*}-0}\lim_{\varepsilon\to0}\left(\frac{|V_{1n}|+|V_{2n}|}{2}-| \mathcal{K}_1- \mathcal{K}_2|\right)
\end{equation*}
in the instant of the confluence and at the point of the contact. Here
$\bar{\theta}=\bar{\theta}(x,y,z,t)$, $(x^*,y^*,z^*)$ is the point of the contact of the free
boundaries, and $\mathcal{K}_1$, $\mathcal{K}_2$ are the principle curvatures of the free boundaries
at the point of the contact.

\begin{figure}
\begin{center}
\includegraphics[width=7cm]{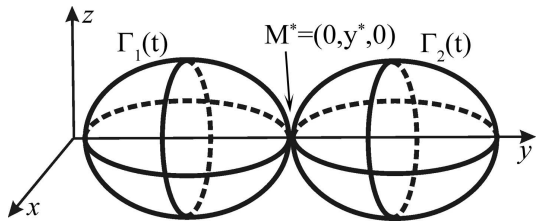}
\end{center}
\caption{Confluence of the free boundaries
in the case of the curvatures of the free boundaries have
the unlike signs.}\label{F2}
\end{figure}

Let us consider another situation. We assume that the solution of our problem is symmetric about
$y$-axe (see Fig.\ref{F2}). By $y_{1}(t)$ and $y_{2}(t)$ we denote the points in which the free
boundaries cross the $y$-axe and these points is situated on the shorts distance one from another.
We assume that $M^*=(0,y^*,0)$ is the point of the contact of the free boundaries. The principle
difference between this situation and the situation in Fig.\ref{F1} is that in the instant of the
contact $t^*$ the principle curvatures of the free boundaries have the unlike signs at the point of
the contact $M^*$. In this case the amplitude of the jump of the temperature is determined by
formula
\begin{equation*}
\left.\left[\bar{\theta}\right]\right|_{(x,y,z)=(0,y^*,0),\ t=t^*} =-\lim_{t\to
t^{*}-0}\lim_{\varepsilon\to 0}
\left(\frac{|V_{1n}|+|V_{2n}|}{2}+|\mathcal{K}_1|+|\mathcal{K}_2|\right).
\end{equation*}

The mentioned soliton type disturbance is observed in the numerical experiments, see
Fig.\ref{F:tminjump},~\ref{F:tt}.

We note that beside the our own papers we know only the single paper of A.M.~Meirmanov and
B.A.~Zaltsman \cite{MZ}. In this paper the problem of the confluence of the free boundaries in the
case of the Hele-Shaw problem is derived and this problem is the special case of the
Stefan-Gibbs--Thomson problem. The regularized (not limit) problem is considered in \cite{MZ} and in
this case the effect of the soliton type disturbance is not shown. The analog of this fact is the
following problem. Let us consider the heat equation in a rectilinear segment and with not zero
($=1$) Dirichlet's initial condition. Clearly, the solution of this problem is zero. However, the
numerical solution is not identical zero for the different scheme with a node in the point $x=0$.

\section{DIFFERENT SCHEME}

The choice of the different scheme for system (\ref{PhHeat}),
(\ref{PhAC}) is based on some ideas that are sufficiently general
for solving nonlinear equations. Namely,
we calculate the heat equation (\ref{PhHeat}) at the next time
step and, consequently split system (\ref{PhHeat}), (\ref{PhAC}).
We associate Eq. (\ref{PhAC}) with implicit different scheme
\begin{align}\label{meshU}
&\varepsilon
u^{k+1}_{\bar{t}}-\varepsilon\sum_{i=1}^{n}u^{k+1}_{r_{i}\bar{r}_{i}}-
2\varepsilon\sum_{i=1}^{n}\frac{u^{k+1}_{r_{i}}}{r_{i}}=\\
&=\frac{1}{\varepsilon} \left[
u^{k+1}-\left((u^{k})^{3}+3(u^{k})^{2}(u^{k+1}-u^{k})\right)\right]
+\varkappa \sum_{i=1}^{n}\frac{\sigma^{k}_{i}}{r_{i}}.\notag
\end{align}
Here $u^{k}$ is the mesh order function at the $k$-th time step,
$\sigma^{k}$ is the mesh "temperature"\ at the $k$-th time step.

The first term (the term in the square  brackets) in the right
hand part of Eq. (\ref{meshU}) is obtained as following. We denote
\begin{equation*}
F(u)=u-u^{3}.
\end{equation*}
We linearize the mesh function $F(u^{k+1})$ as following
\begin{equation*}
F(u^{k+1})\approx F(u^{k})+F'_{u}(u^{k})(u^{k+1}-u^{k}).
\end{equation*}

System (\ref{meshU}) is completed with initial and boundary
conditions in the points $r=R_{j}$, $j=1,2$. It is clear that Eqs.
(\ref{meshU}) is the three-point equations relatively to $u^{k+1}$
and these equations are solved by the sweep method.

To calculate the mesh "temperature"\ $\sigma^{k+1}$ we use the
standard different scheme for the heat equation (\ref{PhHeat}).
Namely,
\begin{equation}\label{meshT}
\sigma_{\bar{t}}^{k+1}-\sum_{i=1}^{n}\sigma^{k+1}_{r_{i}\bar{r}_{i}}=
-\sum_{i=1}^{n}r_{i}\sigma^{k+1}_{i\bar{t}}.
\end{equation}
Equations (\ref{meshT}) is also solved by the sweep method.

We use the main segment $r\in[1,2]$ (i.e., $R_{1}=1$, $R_{2}=2$)
for the numerical simulation of the process of the confluence.

\begin{figure}
\begin{center}
\includegraphics[width=7cm]{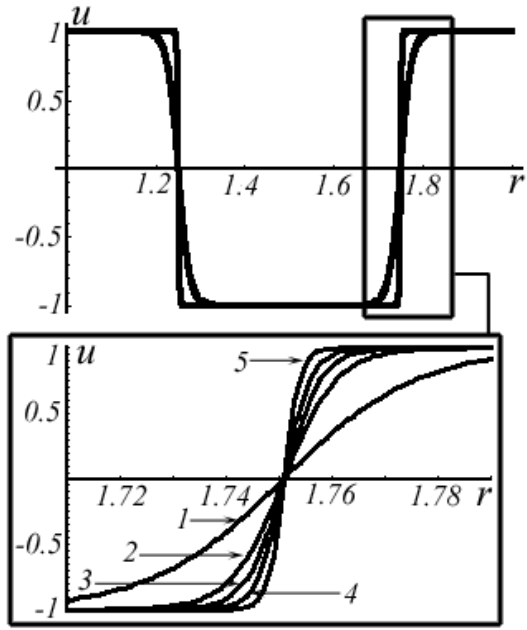}
\end{center}
\caption{The initial data for the order function $u$ for the
different values of $\varepsilon$: 1. $\varepsilon=0.025$, 2.
$\varepsilon=0.01$, 3. $\varepsilon=0.007$, 4.
$\varepsilon=0.005$, 5. $\varepsilon=0.003$.}\label{F:u0}
\end{figure}

\begin{figure}
\begin{center}
\includegraphics[width=7cm]{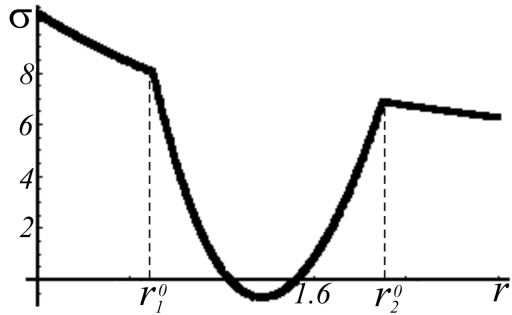}
\end{center}
\caption{The initial data for the temperature $\sigma$. The
initial positions of the free boundaries are $r_{1}^{0}=1.25$,
$r_{2}^{0}=1.75$}\label{F:t0}
\end{figure}

The initial data we choose as provided by the structure of the
weak asymptotic solution (\ref{AS:2}). We use the function
\begin{equation}\label{u0}
u^{0}=1+\tanh\left(\frac{r_{1}^{0}-r}{\varepsilon}\right)+
\tanh\left(\frac{r-r_{2}^{0}}{\varepsilon}\right)
\end{equation}
as the initial data for the order function $u$. Here
$r_{1}^{0}=1.25$, $r_{2}^{0}=1.75$ determine the initial position of the
free boundaries (see Fig.\ref{F:u0},\ref{F:t0}).

The initial data $\sigma^{0}$ for the "temperature"\ is taken in
the form in Fig.\ref{F:t0} as provided by model (\ref{ModT:1}).
Namely, the function $\sigma^{0}$ is parabolic in the domain with
phase "$-$"\ (between the free boundaries), and the function
$\sigma^{0}$ is linear in the domain with phase "$+$". At the same
time, the dependence between the initial positions and the initial
velocities of the free boundaries and the values of $\sigma^{0}$
is determined by formula (\ref{SGT:1})in the boundary points
$r_{1}^{0}$, $r_{2}^{0}$.

It is clear that the considered initial data differ from the exact
solution of the problem. Nevertheless it is known that the
solutions of problem (\ref{PhHeat}), (\ref{PhAC}) converge to
solutions of the limit Stefan-Gibbs-Thomson problem as
$\varepsilon\to 0$ and on the sufficiently general assumptions.
More other, the solution of the Cauchi problem come to
self-similar regime by the widely known behaviours of the
semi-linear parabolic equations (see \cite{OR}). Beside it, the
formulas for the asymptotic solution are the deformations of the
formulas for the semi-similar solution.

\section{THE RESULTS OF NUMERICAL SIMULATION}

Let us consider the results of the numerical simulation of the
process of the confluence of the free boundaries.

Some graphics given below illustrate the common behaviors of the solution of system (\ref{PhHeat}),
(\ref{PhAC}) and confirm the certainty of the numerical results. We simulate for the different
(decreasing) values of the parameter $\varepsilon$, $\varepsilon=0.025$, $\varepsilon=0.01$,
$\varepsilon=0.007$, 4. $\varepsilon=0.005$, $\varepsilon=0.003$. The mesh is equal $h=10^{-3}$  and
the time step is equal $\tau=10^{-5}$.

\begin{figure}
\begin{center}
\includegraphics[width=8cm]{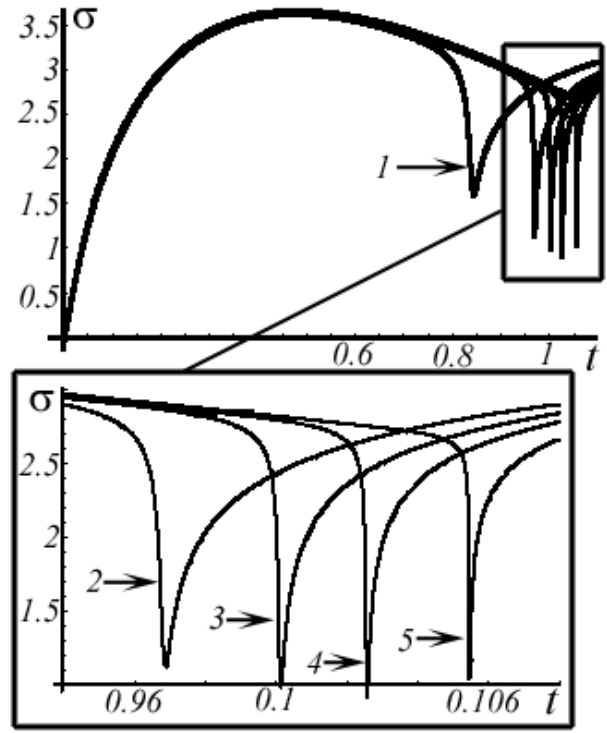}
\end{center}
\caption{The jump of the temperature $\sigma$ in the neighborhood of the point of the contact of the
free boundaries ($r^{*}\approx r=1.42$): 1. $\varepsilon=0.025$, 2. $\varepsilon=0.01$, 3.
$\varepsilon=0.007$, 4. $\varepsilon=0.005$, 5. $\varepsilon=0.003$.}\label{F:tjump}
\end{figure}

\begin{figure}
\begin{center}
\includegraphics[width=8cm]{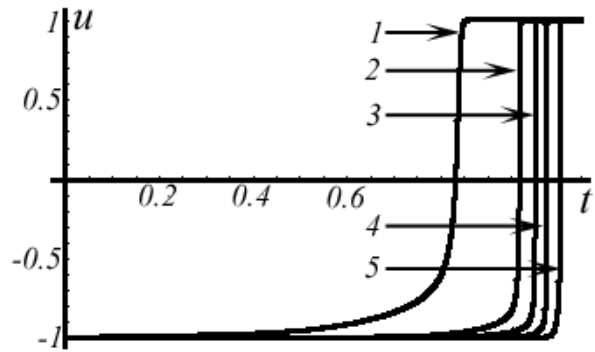}
\end{center}
\caption{The order function $u$ in the neighborhood of the point of the contact of the free
boundaries ($r^{*}\approx 1.42$): 1. $\varepsilon=0.025$, 2. $\varepsilon=0.01$, 3.
$\varepsilon=0.007$, 4. $\varepsilon=0.005$, 5. $\varepsilon=0.003$.}\label{F:ujump}
\end{figure}

In Fig.\ref{F:tjump},~\ref{F:ujump} are shown the profiles of the
temperature $\sigma$ and the order function $u$ in the
neighborhood of the point of the contact of the free boundaries
(for fixed $r^{*}\approx r=1.42$) for different values of the
parameter $\varepsilon$.

\begin{figure}
\begin{center}
\includegraphics[width=7cm]{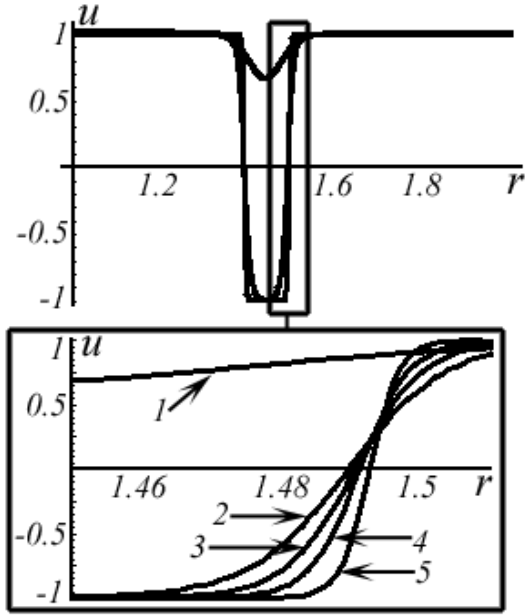}
\end{center}
\caption{The order function $u$ at $t=0.082$: 1.
$\varepsilon=0.025$, 2. $\varepsilon=0.01$, 3.
$\varepsilon=0.007$, 4. $\varepsilon=0.005$, 5.
$\varepsilon=0.003$.}\label{F:ur}
\end{figure}

From the graphics in Fig.\ref{F:tjump} we see that the "width"\ of the neighborhood of the instant
of time of the jump of the temperature decrease (it is proportional to $\varepsilon$) owing to
decreasing $\varepsilon$. At the same time, the jump occurs at the the different instant of time for
the different values of the parameter $\varepsilon$. Namely, in dependence of the decreasing of
$\varepsilon$ the sequence of the instant of time $t_{min}$ (in which the jump becomes minimum)
increases, see also Table~\ref{T:tmin}. Beside it, from Table~\ref{T:tmin} we see that the intervals
grow short between the neighboring instants $t_{min}$ as the parameter $\varepsilon$ decreases. The
cause of this fact is that the width of the transition zone of the order function $u$ is different
for different values of $\varepsilon$, see formula (\ref{u0}). Namely, the width of transition zone
is proportional $\varepsilon$ (as larger the value of $\varepsilon$, as wider the transition zone),
see Fig.~\ref{F:u0} and \cite{Cag1, DOR}. As a result of this fact the contact of the transition
zones (and, consequently the beginning  of the confluence of the free boundaries) occurs previously
for the simulation process with larger $\varepsilon$, see Fig.~\ref{F:ur}. Clearly, the width of the
confluence is also proportional to $\varepsilon$ respect to $r$. From Fig.~\ref{F:ur} we see that
for the fixed instant of time $t=0.082$ the free boundaries corresponding to the graphic 1
($\varepsilon=0.025$) complete the confluence, but at the same time the free boundaries
corresponding another graphics (for smaller $\varepsilon$) are on the sufficiently large distance.

From Fig.~\ref{F:ujump} we see that as smaller $\varepsilon$ as
quickly transit occurs from the phase "$-$"\ to the phase "$+$".
It is clear that $u(r^*,t)\to \mathrm{sign}(t-t_{min})$ as
$\varepsilon\to 0$.

\begin{figure}
\begin{center}
\includegraphics[width=7cm]{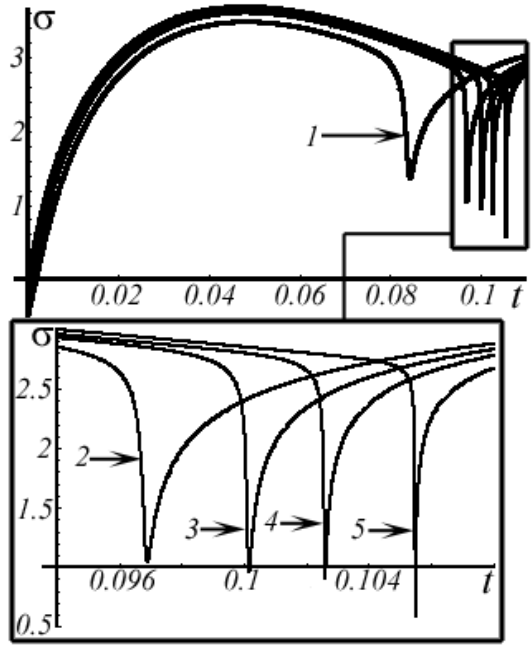}
\end{center}
\caption{The jump of the temperature $\sigma$: 1.
$\varepsilon=0.025$, $r=1.448$, 2. $\varepsilon=0.01$, $r=1.427$,
3. $\varepsilon=0.007$, $r=1.422$, 4. $\varepsilon=0.005$,
$r=1.419$, 5. $\varepsilon=0.003$, $r=1.415$.}\label{F:tminjump}
\end{figure}

\begin{table}
\begin{center}
\begin{tabular}{|c|c|c|c|}
\hline $\varepsilon$ & $t_{min}$ & $r_{min}$ &
$\sigma(r_{min},t_{min},\varepsilon)$ \\ \hline 0.025 & 0.08432 &
1.448 & 1.373699 \\ \hline 0.01 & 0.09690 & 1.427 & 1.051388 \\
\hline 0.007 & 0.10016 & 1.422 & 0.967825 \\ \hline 0.005 &
0.10260 & 1.419 & 0.902497 \\ \hline
0.003 & 0.10549 & 1.415 & 0.587919 \\
\hline
\end{tabular}
\caption{Coordinate and time dependence of the jump of the
temperature $\sigma$. $r_{min}$ and $t_{min}$ are values of the
coordinate $r$ and time $t$ at which the jump becomes
minimum.}\label{T:tmin}
\end{center}
\end{table}

Beside it, from Fig.\ref{F:tjump} we see that the minimum of the
jump for the graphic 5 ($\varepsilon=0.003$) is smaller then the
minimum of the jump for the graphic 4 ($\varepsilon=0.005$). But
in real this is not true and we observe reversed dependence, see
Table~\ref{T:tmin} and Fig.\ref{F:tminjump}. Here are shown the
graphics of the temperature $\sigma$ in the fixed point
$r=r_{min}$, where $r_{min}$ is the value of the coordinate $r$,
in which the jump of the temperature becomes zero. This fact means
that the shift of the minimum of the jump of the temperature
$\sigma$ depend on $\varepsilon$ in the coordinate $r$, see
Table~\ref{T:tmin}. From Table~\ref{T:tmin} we see that this shift
is not lager the difference between the width of the transition
zones (difference between the the corresponding values of the
parameter $\varepsilon$) for the corresponding graphics, see
Fig.\ref{F:tminjump}.

From Table~\ref{T:tmin} we see that the intervals between the neighborhood values of $t_{min}$
decrease as $\varepsilon$ decreases. More other, this intervals is the order of $\varepsilon$. At
the same time the distance between the neighborhood points $r_{min}$ grows short as $\varepsilon$
decreases.

\begin{figure}
\begin{center}
\includegraphics[width=8cm]{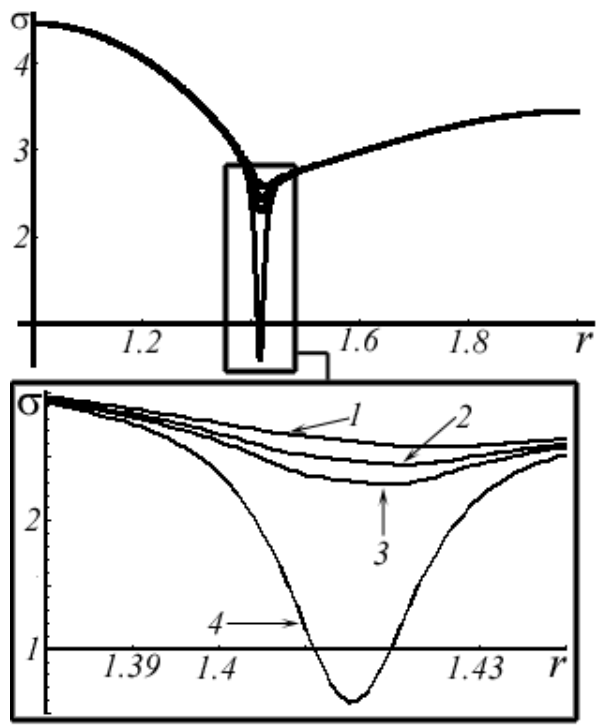}
\end{center}
\caption{The dynamic of the jump of the temperature $\sigma$ for
$\varepsilon = 0.003$: 1. $t=0.1052$, 2. $t=0.10539$, 3.
$t=0.10544$, 4. $t=t_{min}=0.10549$.}\label{F:tt}
\end{figure}

The dynamic of the jump of the temperature at time is shown in
Fig.~\ref{F:tt} for $\varepsilon = 0.003$.

So we can conclude that the process of the confluence of the free
boundaries is localized in the coordinate and at the time.
Consequently, the soliton type disturbance of the temperature is
localized in the coordinate and at the time (the width of the
localization is proportional  to$\varepsilon$), see Fig.~\ref{F0}.

We denote that the values of the parameter $\varepsilon>0.003$ is not possible in the our
simulation. For $\varepsilon=0.003$ and $h=0.001$ we obtain that the three nodes (the minimal number
of the nodes that necessary to calculate the second difference derivation, see) of the mesh appear
in the transition zone, see (\ref{meshU}), (\ref{meshT}). The computer simulation shows that the
numerical solution is unstable in the neighborhood of the confluence of the free boundaries for
$\varepsilon=0.002$ and the different scheme do not give solution for $\varepsilon=0.001$.

\begin{figure}
\begin{center}
\includegraphics[width=8cm]{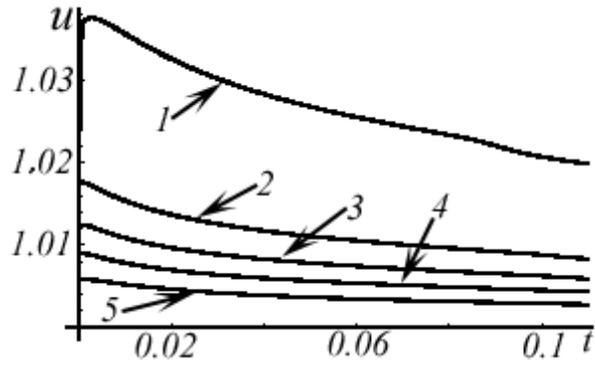}
\end{center}
\caption{The order function $u$ in $r=1.15$: 1.
$\varepsilon=0.025$, 2. $\varepsilon=0.01$, 3.
$\varepsilon=0.007$, 4. $\varepsilon=0.005$, 5.
$\varepsilon=0.003$.}\label{F:u1.15}
\end{figure}

\begin{figure}
\begin{center}
\includegraphics[width=7cm]{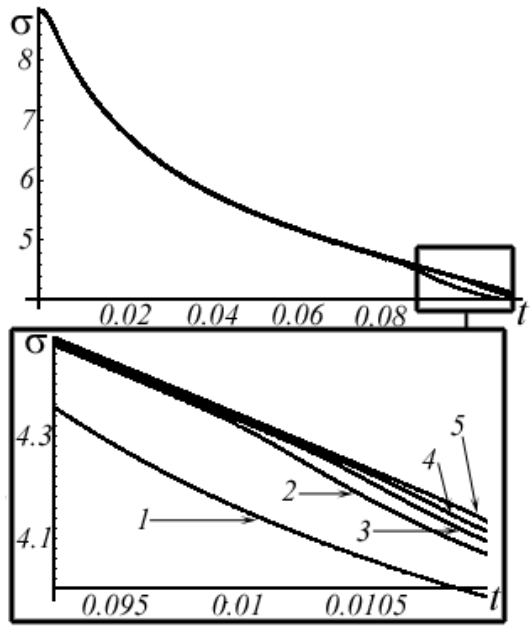}
\end{center}
\caption{The temperature $\sigma$ in $r=1.15$: 1.
$\varepsilon=0.025$, 2. $\varepsilon=0.01$, 3.
$\varepsilon=0.007$, 4. $\varepsilon=0.005$, 5.
$\varepsilon=0.003$.}\label{F:t1.15}
\end{figure}

The series of the graphics demonstrates the stability of the numerical solution of system
(\ref{PhHeat}), (\ref{PhAC}) respect to $\varepsilon$. In Fig.\ref{F:u1.15},~\ref{F:t1.15} are shown
the graphics of the dependence of the order function $u$ and of the temperature $\sigma$ on time in
the fixed point $r=1.15$. We note that the analogous results are correct for the another value of
$r$ except the point of the contact of the free boundary.

From the graphics in Fig.\ref{F:u1.15} we see that at the beginning ($t\leqslant 0.01$) the solution
$u$ undergo the sharp jump and thereafter stables (undergo to near semi-similar regime). At the same
time, this stabilization is wavelike. More other, the initial interval of the time in which the
solution undergo the most strong change (jump) grows short as $\varepsilon$ decreases. More other
the distance between the neighborhood graphics grows short as $\varepsilon$ decreases and this
distance is proportional.

From the graphics in Fig.\ref{F:t1.15} we see that the temperature
$\sigma$ is stable respect to the variation of $\varepsilon$. For
$t>0.09$ the difference between the graphics  due to the process
of the confluence for the graphics for the larger $\varepsilon$
begins early than for the graphics for the smaller~$\varepsilon$.

\begin{figure}
\begin{center}
\includegraphics[width=8cm]{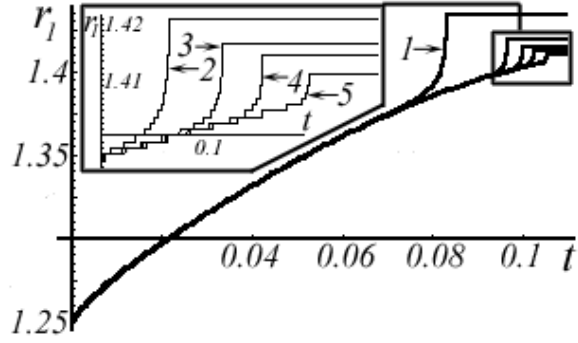}
\end{center}
\caption{The dynamic of the free boundary $r_{1}(t)$ for the
different values of the parameter $\varepsilon$: 1.
$\varepsilon=0.025$, 2. $\varepsilon=0.01$, 3.
$\varepsilon=0.007$, 4. $\varepsilon=0.005$, 5.
$\varepsilon=0.003$.}\label{F:r1}
\end{figure}

\begin{figure}
\begin{center}
\includegraphics[width=8cm]{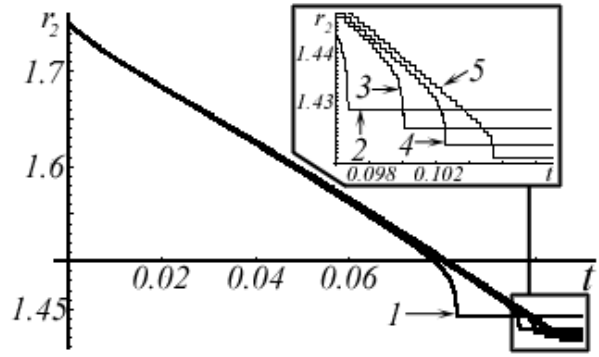}
\end{center}
\caption{The dynamic of the free boundary $r_{2}(t)$ for the
different values of the parameter $\varepsilon$: 1.
$\varepsilon=0.025$, 2. $\varepsilon=0.01$, 3.
$\varepsilon=0.007$, 4. $\varepsilon=0.005$, 5.
$\varepsilon=0.003$.}\label{F:r2}
\end{figure}

The last series of the numerical experiments deals with the
verification of the analytical formula for the amplitude of the
jump of the temperature. Namely, in paper~\cite{DR} is obtained
the formula for the amplitude of the jump
\begin{equation*}
\left.\left[\bar{\theta}\right]\right|_{r=r^*,\ t=t^*}
=-\frac{r_{10}'-r_{20}'}{2},
\end{equation*}
or taken into account (\ref{sigma}) we have
\begin{equation}\label{sjump}
\left.\left[\bar{\sigma}\right]\right|_{r=r^*,\ t=t^*}
=-\lim_{t\to t^{*}-0}\lim_{\varepsilon\to
0}\frac{2r_{1}(t,\varepsilon)r_{1}'(t,\varepsilon)-2r_{2}(t,\varepsilon)r_{2}'(t,\varepsilon)}{4}.
\end{equation}

In Fig.\ref{F:r1},~\ref{F:r2} the dependence of the free boundaries $r_{1}(t)$ and $r_{2}(t)$ on
time is shown for the different values of the parameter $\varepsilon$. We see that the functions
$r_{1}(t)$ and $r_{2}(t)$ is linear in the initial interval of the shifting of the free boundaries.
The nonlinearity of the confluence depends on the shifting of the free boundary near the point of
the contact and the lines are distorted. We see that as smaller $\varepsilon$ as smaller the
neighborhood of the instant $t^{*}$ in which the distortion is sensed.

The functions $r_{i0}(t)$, $i=1,2$ are determined as the limits
\begin{equation*}
r_{i0}=\lim_{\varepsilon\to 0}r_{i}(t,\varepsilon)
\end{equation*}
in the asymptotic formulas. So in the capacity of
$\frac{d}{dt}r_{i0}(t^{*})$ we should take the constants which are
obtained by differentiation of the functions
$r_{i}(t,\varepsilon)$ in the interval where this functions the
most close to the line.

We note that the positions of the free boundaries is fuzzy respect to $\varepsilon$ as
$\varepsilon>0$ (the width of the transition zone is proportional to $\varepsilon$). Therefore to
disclose the behavior of the shifting the free boundaries we consider the dynamic of the point of
the crossing the order function $u(r,t)$ and the axis $r$, i.e. we trace the dynamic of the points
$r(t)$ in which $u(r(t),t)=0$ for any fixed $t$, see Fig.\ref{F:r1},~\ref{F:r2}. The instants of
time $t^*$ in which the graphics of $r_{1}(t)$ and $r_{2}(t)$ become the horizontal lines correspond
to the instants of time in which the order function $u(r,t)$ becomes positive. This instants of time
naturally do not equal to instants of time $t_{min}$ in which the amplitude of the jump of the
temperature is maximal, see Table~\ref{T:tmin}.

To calculate the amplitude of the jump of the temperature we
determine the minimal value of the temperature $\sigma_{min}$
corresponded to the value of the coordinate $r_{min}$, see
Table~\ref{T:tmin}. For example in Fig.~\ref{F:tt}, this value of
the coordinate is become on the  curve 4. The amplitude is equal
to the absolute value of the difference between $\sigma_{min}$ and
the value of the temperature in the point $r_{min}$ and at the
instant of time corresponded to the start of the fast changing of
the temperature in this point (in Fig.~\ref{F:tt} the curve 3
gives this value of the temperature).

The results of the verification of the formula (\ref{sjump}) are
showed in Table~\ref{T:tjump}.

\begin{table}
\begin{center}
\begin{tabular}{|c|c|c|c|c|c|c|}
\hline $\varepsilon$ & $r_{1}$ & $r'_{1}$ & $r_{2}$ & $r'_{2}$ &
$\left[\sigma_{an}\right]$ & $\left[\sigma_{cal}\right]$ \\
\hline 0.025 & 1.394 & 1.4 & 1.485 & -3.2 &
3.01448 & 1.5 \\
\hline 0.01 & 1.402 & 1 & 1.446 & -3.2 &
3.3518 & 1.7 \\
\hline 0.007 & 1.403 & 0.8 & 1.438 & -3.2 &
2.86 & 1.9 \\
\hline 0.005 & 1.403 & 0.6 & 1.430 & -3.2 &
2.7 & 2.1 \\
\hline 0.003 & 1.405 & 0.6 & 1.422 & -3.2 &
2.69 & 2.7 \\
\hline
\end{tabular}
\caption{The value of the amplitude of the jump of the temperature
$\sigma$ obtained by analytical formula
($\left[\sigma_{an}\right]$) and by numerical simulation
($\left[\sigma_{cal}\right]$).}\label{T:tjump}
\end{center}
\end{table}
In column $\left[\sigma_{an}\right]$ the data for the jump are
given obtained by formula (\ref{sigma}), and in column
$\left[\sigma_{cal}\right]$ the data for the jump are given
calculated by numerical simulation. From Table~\ref{T:tjump} we
see that the analytical data converge to the numerical data as
$\varepsilon$ decreases.


\begin{thebibliography}{99}

\bibitem{Cag1}
\textsc{G.Caginalp} \emph{An analysis of a phase field model of a
free boundary.} Arch. Rat. Mech. Anal. 92 (1986), 205--245.

\bibitem{DOR}
\textsc{V.G.Danilov, G.A.Omel'yanov \& E.V.Radkevich}
\emph{Asymptotic solution of a phase field system and the modified
Stefan problem.} 1995, Differential'nie Uravneniya 31(3), 483--491
(in Russian).(English translation in Differential Equations 31(3),
1993.

\bibitem{D}
\textsc{V.G.Danilov} \emph{Weak asymptotic solution of phase-field
system in the case of confluence of free boundaries in Stefan
problem with undercooling.} Euro. Journ. Appl. Math. (2007) v.18,
pp 537-569.

\bibitem{DR}
\textsc{V.G.Danilov, V.Yu.Rudnev} \emph{Confluence of the nonlinear waves in the
Stefan--Gibbs-Thomson problem.} Free boundary problems. To be published.

\bibitem{OR}
\textsc{G.A.Omel'yanov and V.Yu.Rudnev} \emph{Interaction of free
boundaries  in the modified Stefan problem.} Nonlinear Phenomena
in Complex Systems, 7:3 (2004) 227 - 237.

\bibitem{PS}
\textsc{P.I.~Plotnikov and V.N.~Starovoitov} \emph{Stefan problem
as the limit of the phase field system.}, Differential Equations
29 (1993), 461--471.

\bibitem{MZ}
\textsc{A.Meirmanov, B.Zaltzman} \emph{Global in time solution to
the Hele-Shaw problem with a change of topology.} EJAM, Vol. 13,
pp. 431-447, 2002.

\bibitem{M}
\textsc{A.Meirmanov} \emph{The Stefan Problem.} Berlin . New York:
Walter de Gruyter, 1992.


\end{thebibliography}
\end{document}